\documentclass[12pt]{amsart}
\usepackage{amscd,amsmath}
\newcommand{\la}{\lambda}
\newcommand{\al}{\alpha}
\newcommand{\be}{\beta}
\newcommand{\ga}{\gamma}
\newcommand{\dd}{\delta}
\newcommand{\e}{\hat\eta} 
\newcommand{\ka}{K\"ahler}

\newcommand{\ov}{\overline}

\newcommand{\E}{\mathcal{E}}
\newcommand{\F}{\mathcal{F}}

\newcommand{\f}{\varphi}
\newcommand{\vv}{\overrightarrow}   

\newcommand{\CC}{\mathbb{C}}   
\newcommand{\HH}{\mathbb{H}}   
\newcommand{\RR}{\mathbb{R}}   
\newcommand{\ZZ}{\mathbb{Z}}

\newcommand{\CP}[1]{\mathbb{C}P^{#1}}
\newcommand{\HP}[1]{\mathbb{H}P^{#1}}

\numberwithin{equation}{section}

\newtheorem{te}{Theorem}[section]
\newtheorem{pr}{Proposition}[section]
\newtheorem{co}{Corollary}[section]

\theoremstyle{definition}
    
\theoremstyle{remark}
\newtheorem{re}{Remark}[section]

\begin{document}

\title[ Moment maps and induced Hopf bundles]
{On some moment maps and induced Hopf
bundles in the quaternionic projective space}
\author{Liviu Ornea}
\address{University of Bucharest\\ Faculty of Mathematics\\ 14 Academiei
str.\\ 70109
Bucharest, Romania}
\email{lornea@imar.ro}
\author{Paolo Piccinni}
\address{       Universit\`{a} degli Studi di Roma
''La Sapienza''\\ Piazzale Aldo Moro 2\\ 
I-00185 Roma, Italia  }
\email{piccinni@mat.uniroma1.it}


\maketitle

\begin{abstract}  We describe a diagram containing the zero sets of the
moment maps associated to the diagonal
$U(1)$ and $Sp(1)$ actions on the quaternionic projective space $\HP{n}$.
These sets are related both to focal sets of
submanifolds and to Sasakian-Einstein structures on induced Hopf bundles.
As an application, we construct a 
complex structure on the Stiefel manifolds
$V_2({\CC}^{n+1})$ and  $\widetilde V_4({\RR}^{n+1})$, the one on the
former manifold not being compatible
with its known hypercomplex structure. 
\end{abstract}

\section{Introduction}

The main motivation for the present work is to understand some aspects of
the Riemannian geometry of the focal set $Foc_{\HP{n}}\CP{n}$, i.e. of the
set of points in the
quaternionic projective space $\HP{n}$ that are critical values of the
normal exponential map with
respect to the totally geodesic submanifold $\CP{n} \subset \HP{n}$. Our
starting point is
the Habilitationsschrift of J. Berndt \cite{Ber 1}, \cite{Ber 2}, where in
particular is proved 
that $Foc_{\HP{n}}\CP{n}$ fibers in circles over $Gr_2(\CC^{n+1})$, the Grassmannian
of 2-planes in
$\CC^{n+1}$. In the simplest cases $n=1,2,$ these focal points distribute on spheres
$S^1$, $S^5,$ respectively
and $n=3$ thus seems to be the first significant case. Indeed, since
$Foc_{\HP{3}}\CP{3}$ is a circle bundle over the Klein quadric
$Gr_2(\CC^4)\cong \CC Q^4\subset \CP{5}$, both the classical projective geometry of
$\CC Q^4$ (described for example in \cite{WW}, pp. 26-35) and its two quaternion
K\"ahler structures \cite{Mor} can be related to the total space
$Foc_{\HP{3}}\CP{3}$. More generally, the twofold K\"ahler-Einstein and quaternion
K\"ahler properties of the Grassmannian $Gr_2(\CC^{n+1})$ suggested to us a link
between the focal
set $Foc_{\HP{n}}\CP{n}$ and the Sasakian-Einstein geometry of the induced Hopf
bundles considered in our previous work \cite{Or-Pi 3}. 

This paper begins by giving an explicit identification of the focal set 
$Foc_{\HP{n}}\CP{n}$ with the zero level set of the moment map
$\mu$ associated to the diagonal $U(1)$ - action on $\HP{n}$.
This is the simplest isometric action giving examples of reductions in quaternion
K\"ahler
geometry, and it is well-known that the reduced manifold is precisely the quaternion
K\"ahler Wolf space $Gr_2(\CC^{n+1})$ \cite{Ga},
\cite{Ga-La}. According to the notation in use, we write $\HP{n}/\! /\! /U(1) =
Gr_2(\CC^{n+1})$ to
indicate this reduction procedure; this has been proved to be the unique one with
respect to isometric circle
actions and positive quaternion K\"ahler manifolds \cite{Bat 3},
\cite{Bat 4}. Next, the diagonal $Sp(1)$-action on $\HP{n}$ has as reduction the
quaternion K\"ahler manifold $\HP{n}/\! /\! /Sp(1) = \widetilde Gr_4(\RR^{n+1})$,
the Grassmannian of oriented 4-planes
in $\RR^{n+1}$. The zero set of the corresponding moment map $\nu: \HP{n}
\rightarrow \HH$ can also
be related to focal sets in $\HP{n}$. (For these identifications cf. Theorems 3.1,
3.2).

On the other hand, the zero level sets $\mu^{-1}(0)$ and $\nu^{-1}(0)$ of these
moment maps $\mu$ and $\nu$ can be
identified  with the total spaces of some induced Hopf $S^1$-bundles. This enables
us to define on them a
Sasakian-Einstein and a 3-Sasakian metric, respectively, and hence to find a diagram
where these zero sets $\mu^{-1}(0)$ and
$\nu^{-1}(0)$ are base spaces of $S^3$-bundles, now induced by the Hopf fibration
$S^{4n+3} \rightarrow \HP{n}$. The corresponding total spaces are the Stiefel
manifolds
$V_2(\CC^{n+1})$ and $\widetilde V_4(\RR^{n+1})$, respectively. The $3$-Sasakian
structure of the fibers $S^3$ and the
Sasakian structure of the base spaces
$\mu^{-1}(0)$ and
$\nu^{-1}(0)$ allow to construct a complex structure on $V_2(\CC^{n+1})$ and on
$\widetilde V_4(\RR^{n+1})$, whose definition is very much in the Calabi-Eckmann spirit. On the former Stiefel manifold, this
complex structure is not compatible with its standard hypercomplex structure, obtained in \cite{Bat 2}, \cite{Bo-Ga-Ma 2},
\cite{Jo}. On the other hand, both $V_2(\CC^{n+1})$ and  $\widetilde V_4(\RR^{n+1})$ are total spaces of framed bundles in
Hopf surfaces over K\"ahler Einstein manifolds. Then these complex structures can be seen to belong to one-parameter
families, also suitable for some exceptional cases, as shown in a forthcoming paper, \cite{Or-Pi 4}. 

{\bf Acknowledgements.} This paper was written while the first named author was
visiting the University La Sapienza of Rome 
and the Max Planck Institut f\"ur Mathematik in Bonn. He thanks  both institutions
for financial support and for 
hospitality. The second author acknowledges support from MURST project ''Propriet\`a
geometriche delle variet\`a reali e complesse''.
\par Both authors thank Victor Vuletescu for some helpful discussions, and John C.
Wood who accepted to correct our English and patiently read the whole paper. 
\bigskip

\section{Preliminaries} 

We collect here some definitions and basic facts that will be used throughout the paper. 
\par Let $(M,g)$ be a Riemannian manifold and let  $N$ be an isometrically immersed
submanifold. The critical values of the restriction  $exp_{T^\perp N}$  of the
exponential map of
$M$ to  the normal bundle of $N$ are called 
\emph{ focal points of $N$} and the set of all focal points, here denoted by
$Foc_MN$, is called
\emph{the focal set of $N$} (see for example \cite{CE} p. 23, \cite{doC}, p. 227,
\cite{ON} p. 283).
The focal set of a submanifold may not be a submanifold: for smooth plane
curves and for regular
surfaces in $\RR^3$ the focal sets are respectively the evolutes and the surfaces of
the centres, and both
of them can have singular points, cf. \cite{doC}, pp. 237-238 and p. 232. However, for the  
examples considered in the present paper all the focal sets turn out to be smooth. Indeed, it seems that
not many focal sets of Riemannian submanifolds have been explicitly determined, but
we can mention totally geodesic spheres in spheres 
\cite{Ot}, and hypersurfaces in space forms \cite{Ce-Ry 1}, \cite{Ce-Ry 2}. 

We are mainly interested in $Foc_{\HP{n}}\CP{n}$, the focal set of
the totally geodesic
$\CP{n}$ in the quaternionic projective space $\HP{n}$. This focal set appears in  
J. Berndt's work \cite{Ber 1}, \cite{Ber 2}, denoted there by $Q^n$, and studied in
connection with both the complex K\"ahler and the quaternion \ka\ structure of
$Gr_2(\CC^{n+1})$, the Grassmannian of complex two-planes in $\CC^{n+1}$. 

Two geometric structures which appear naturally in our context are the Sasakian and
the $3$-Sasakian ones. We briefly recall the definitions, referring the reader to
the survey \cite{Bo-Ga 2} for further information. 
\par A {\it Sasakian} manifold is a $(2n+1)$-dimensional Riemannian manifold $(S,g)$
equipped with a unitary Killing vector field $\xi$ such that the field of
endomorphisms
$\f:=\nabla\xi$\ satisfies the differential equation
$$\nabla\f=Id\otimes\eta-g\otimes\xi,$$ where
$\nabla$ is the Levi-Civita connection of $g$ and $\eta$ is the dual 1-form of
$\xi$. The data of a
Sasakian structure on the manifold $S$ is easily seen to be equivalent to the
requirement that the cone metric $dr^2+r^2g$, on
$\RR_+\times S$ have holonomy contained in
$U(n+1)$. Note that
$\eta$ is a contact form on $S$, hence $\xi$ is its Reeb field. 
In the simplest example, the Euclidean sphere
$S^{2n+1}$, the Killing vector field is $\xi=-JU$, $J$ being the standard complex
structure of $\CC^{n+1}$ and $U$ the unit outward normal to the sphere. 

More generally, we look at the induced 
Hopf $S^1$-bundle $\pi:V\rightarrow M$, over a smooth submanifold of $\CP{N}$. Its
total space $V$ carries a Sasakian structure
$(V,g,\xi)$ induced from the one of $(S^{2N+1}, can)$. If the Fubini-Study metric of
$\CP{N}$ induces an Einstein metric $h$ with Einstein constant $\al$ on $M$, it can
be seen (cf. \cite{Or-Pi 3}, lemma 1) that the Ricci tensor of the metric
$g$ has the form 
$$Ric(g)=\la g+\mu\eta\otimes\eta,$$ with $\la=\al-2,$ $\mu=\dim V+1-\al$. This
is known as the {\it $\eta$-Einstein} property in Sasakian geometry and, 
following S. Tanno \cite{Ta}, an
$\eta$-Einstein Sasakian metric can be deformed to a Sasakian-Einstein one by setting
\begin{equation}\label{mom1}
g'=Ag+A(A-1)\eta\otimes\eta,
\end{equation}
so that $\xi'= A^{-1}\xi$ with $A=\frac{\la+2}{\dim V+1}$. Now $\pi$ is a Riemannian
submersion with respect to the metrics $h$ and $g$; thus to have a Riemannian
submersion with respect to $g'$, we must consider the scaled metric
$h'=Ah$ on the base $M$.

As Sasakian geometry may be viewed as the odd-dimensional counterpart of \ka\
geometry, the
odd-dimensional counterparts of hyper\ka\ manifolds are $3$-Sasakian ones (cf.
\cite{Bo-Ga 2}). More precisely, a $(4n+3)$-dimensional Riemannian manifold $(S,g)$
is said to be
$3$-Sasakian if it is endowed with three mutually orthogonal unit Killing vector
fields $\xi_1$, $\xi_2$, $\xi_3$, each one defining a Sasakian structure and
satisfying the conditions: $$[\xi_1,\xi_2]=2\xi_3,
\hspace{0.5cm} [\xi_2,\xi_3]=2\xi_1,\hspace{0.5cm}
[\xi_3,\xi_1]=2\xi_2.$$ As above,
an equivalent definition requires that the cone
metric $dr^2+r^2g$ on $\RR_+\times S$ be hyperk\"ahler, i.e. its holonomy group
be contained in
$Sp(n+1)$.
$3$-Sasakian manifolds are necessarily Einstein with positive scalar curvature and
their Einstein constant is $4n+2$. Given a positive quaternion \ka\ manifold $P$,
one constructs its \ka -Einstein twistor space
$Z_P$ and then an $S^1$ bundle over it whose Chern class is, up to torsion, that of
an induced Hopf bundle. The total space
$S$ thus obtained is an 
$SO(3)$-principal bundle over $P$ with $3$-Sasakian structure. 
Moreover, all three fibrations involved are Riemannian submersions.
  
We recall now two basic moment maps of quaternion K\"ahler geometry. Let
$[h_0:h_1:...:h_n]$ be the homogeneous coordinates on
$\HP{n}$: for each $a=0,1,...,n$, we shall write the complex and real
components of $h_a$ as follows: 
\begin{equation}\label{mom2}
h_a = z_a + w_a j=\al_a+\be_ai+\ga_aj+\dd_ak 
\end{equation}
where $z_a=\al_a+\be_ai,$ and $w_a=\ga_a+\dd_ai.$
The first moment map, induced by the diagonal action of $U(1)$ on $\HP{n}$ is at the
hyperk\"ahler level of the bundle $\HH^{n+1}\backslash \{0\}\rightarrow \HP{n}$,
given by: 
\begin{equation}\label{mom3}
\mu:\HH^{n+1}\backslash \{0\}\rightarrow \HH , \quad \mu=\sum_{a=0}^n\ov{h}_aih_a.
\end{equation}
On $\HP{n}$, one has to regard the corresponding moment map $\mu_{qK}$ as a section
of $S^2 H$, the rank $3$ vector bundle of the compatible almost complex structures.
In fact, by using the metric, 
$\mu_{qK}$ appears as a 2-form solution of $$d\mu_{qK}=i_X \Omega,$$ where $\Omega$
is the K\"ahler 4-form and $X$ the Killing vector field generating the $U(1)$-action:
\cite{Ga},
\cite{Bat 1}. However, since the zero level sets $\mu^{-1}(0) \subset \HH^{n+1}$ and
$\mu_{qK}^{-1}(0) \subset
\HP{n}$ correspond to each other in the bundle projection $\HH^{n+1} \backslash 
\{0\} \rightarrow \HP{n}$, we shall refer to the definition
$\mu=\sum_{a=0}^n\ov{h}_aih_a$, whose zero set makes sense also when the $h_a$ are
the homogeneous coordinates of $\HP{n}$ \cite{Bo-Ga 2}. The reduced manifold
$\mu_{qK}^{-1}(0)/U(1)$ turns out to be the quaternion K\"ahler Wolf space
$\frac{SU(n+1)}{S(U(n-1)\times U(2))}\cong Gr_2(\CC^{n+1})$ \cite{Ga}, \cite{Ga-La}. 

The second moment map to be considered, generated by the action of $Sp(1)$, is
(again at the
hyperk\"ahler level): 
\begin{equation}\label{mom4}
\nu:\HH^{n+1}\backslash \{0\}\rightarrow \HH^3 , \quad
\nu=(\sum_{a=0}^n\ov{h}_aih_a, \sum_{a=0}^n\ov{h}_ajh_a,
\sum_{a=0}^n\ov{h}_akh_a), 
\end{equation}
and its corresponding quaternion K\"ahler moment map $\nu_{qK}$ can be viewed as a
triple of 2-forms associated to the frame of Killing vector fields defining the
$Sp(1)$-action. The corresponding reduced manifold $\nu_{qK}^{-1}(0)/Sp(1)$ is now
the Wolf space
$\frac{SO(n+1)}{SO(n-3)\times SO(4)}\cong\widetilde{Gr}_4(\RR^{n+1})$, the
Grassmannian of oriented four-planes of $\RR^{n+1}$.  

\bigskip

\section{Statement of results}

\begin{te}\label{te1}
\par {\rm (i)} The focal set $Foc_{\HP{n}}\CP{n}$ coincides with
$\mu^{-1}(0)$, and it is isometric to 
the total space of the induced Hopf $S^1$-bundle via the
Pl\"ucker embedding $Gr_2(\CC^{n+1}) \hookrightarrow  \CC P^N$. 
\par {\rm (ii)} The metric
$g_1$, induced on $\mu^{-1}(0)$ by the Pl\"ucker embedding allows us to define on
$\mu^{-1}(0)$ a Sasakian Einstein metric
$g$. 
\par {\rm (iii)} The Stiefel manifold $V_2(\CC^{n+1})$ of orthonormal 2-frames in
$\CC^{n+1}$, diffeomorphic to the
total space of the induced Hopf $S^3$-bundle via the embedding 
$\mu^{-1}(0) \subset \HP{n}$, admits an (integrable) complex structure 
$J$, not compatible with
the  standard hypercomplex structure of $V_2(\CC^{n+1})$.
\end{te}

We have a similar statement regarding the moment map $\nu$ induced by the action of
$Sp(1)$ on
$\HP{n}$. We need to consider the following mutually congruent, totally
geodesic embeddings of  $\CP{n}$ in $\HP{n}$.
\begin{equation}\label{cc}
\begin{split}
\CP{n}_i=&\{h\in \HP{n}\; ;\; \ga_a=\dd_a=0, a=1,...,n\},\\
\CP{n}_j=&\{h\in \HP{n}\; ;\; \be_a=\dd_a=0, a=1,...,n\},\\
\CP{n}_k=&\{h\in \HP{n}\; ;\; \be_a=\ga_a=0, a=1,...,n\}:
\end{split}
\end{equation}
here
$\CP{n}_i$ is the standard $\CP{n}$ appearing in Theorem 3.1.

\begin{te}\label{te2}
\par {\rm (i)} The zero level set $\nu^{-1}(0)$ coincides with the intersection
$M=Foc_{\HP{n}}\CP{n}_i\cap Foc_{\HP{n}}\CP{n}_j\cap Foc_{\HP{n}}\CP{n}_k$ and is
isometric to the total space of the induced Hopf $S^1$-bundle over the Fano manifold
$Z_{\widetilde{Gr}_4(\RR^{n+1})}$, via the embeddings
$Z_{\widetilde{Gr}_4(\RR^{n+1})}\hookrightarrow Gr_2(\CC^{n+1})\hookrightarrow 
\CC P^N$, the first of which is defined by regarding
$Z_{\widetilde{Gr}_4(\RR^{n+1})}$ as the space of totally isotropic two-planes in
$\CC^{n+1}$.  This isometry allows the construction of a Sasakian Einstein metric on
$\nu^{-1}(0)$, and identifies it with the homogeneous 3-Sasakian manifold
$SO(n+1)/(SO(n-3)
\times Sp(1))$.
\par {\rm (ii)} The Stiefel manifold
$\widetilde V_4(\RR^{n+1})$ admits an (integrable) complex structure, 
projecting to the complex structure of
$Z_{\widetilde{Gr}_4(\RR^{n+1})}$.
\end{te} 

Here $Z_{\widetilde{Gr}_4(\RR^{n+1})}$ is
the twistor space of the quaternion K\"ahler Wolf space
$\widetilde{Gr}_4(\RR^{n+1})$ given by the $Sp(1)$ reduction on $\HP{n}$.
$Z_{\widetilde{Gr}_4(\RR^{n+1})}$ is known to be a complex submanifold of the
K\"ahler-Einstein Grassmannian $Gr_2(\CC^{n+1})$, see for example
\cite{Kobak} p. 14 or \cite{He-Sa} p. 702.

Statements \ref{te1}, \ref{te2} give, in particular, some fibrations that can
be collected into a diagram as follows.  Here
$V_k(\mathbb{C}^{n+1})$ and $\widetilde V_k(\mathbb{R}^{n+1})$ denote the Stiefel
manifolds of
$k$-frames in $\CC^{n+1}$ and of oriented $k$-frames in $\RR^{n+1}$, respectively.

\begin{co}\label{te3}
There is a commutative diagram 
$$\begin{array}{cccccl}
\widetilde{V}_4(\mathbb{R}^{n+1})&\;  \hookrightarrow\;&
V_2(\mathbb{C}^{n+1})&\; \hookrightarrow\;&S^{4n+3}&\\[1mm]
\Big\downarrow  {S^3}& &\Big\downarrow  {S^3}& &\Big\downarrow  {S^3} \\[2mm]
\nu^{-1}(0)&\; \hookrightarrow\;&\mu^{-1}(0)&
\hookrightarrow\;&\mathbb{H}P^n&\;\;\;\; S^{2N+1}\\[1mm]
\Big\downarrow  {S^1}& &\Big\downarrow  {S^1}& &&\!\!\!\!\!\swarrow \\[2mm]
Z_{\widetilde{Gr}_4(\mathbb{R}^{n+1})}&\;
\hookrightarrow\;&Gr_2(\mathbb{C}^{n+1})&\;
\hookrightarrow\;&\mathbb{C}P^N&\\[2mm]
\Big\downarrow  {S^2}&&&&\;\\[2mm]
\widetilde{Gr}_4(\mathbb{R}^{n+1})&&&&&
\end{array}$$
\noindent of principal $S^1$ and $S^3$-bundles and Riemannian submersions. The zero
level sets 
$\mu^{-1}(0)$ and $\nu^{-1}(0)$ are thus total spaces of induced Hopf $S^1$-bundles
and, as such, are minimal submanifolds of the sphere
$S^{2N+1}$.   
\end{co}
\bigskip

\section{The proofs}

\subsection*{Proof of Theorem 3.1 (i) and (ii)}

The first observation is:
\bigskip
\par \noindent
{\it $Foc_{\HP{1}}\CP{1}$ is the level set of the moment map associated 
to the standard $U(1)$-action on $\HP{1}$.}
\bigskip
\par 
In the homogeneous coordinates $[h_0:h_1]$ of $\HP{1}$, using the notations from
formula (2.1),   
$\CP{1}$ is given by:  
\begin{equation}\label{doi}w_0=0, w_1=0\; ~\text{or}~ \; 
\ga_0=\dd_0=\ga_1=\dd_1=0.
\end{equation}
The pair $\CP{1} \subset \HP{1}$ can be identified with $S^2 \subset S^4$, and the
focal
set of a totally geodesic sphere $S^p$ in $S^n$ is the unit $S^{n-p-1}$ in the
$(n-p)$-dimensional
orthogonal complement of the $\RR^{p+1}$ containing $S^p$
(\cite{Ot}, p. 286).  This identification can be made explicit by using the
coordinate $h=\al+\be i+\ga j+\dd k=h_0h_1^{-1}$ in the affine line $h_1\neq 0$,
whose corresponding real coordinates are: 
\begin{equation*}
\begin{split}
\al&=\frac{\al_0\al_1+\be_0\be_1+\ga_0\ga_1+\dd_0\dd_1}{r^2},\; 
\be=\frac{\al_1\be_0-\al_0\be_1-\ga_0\dd_1+\ga_1\dd_0}{r^2},\\
\ga&=\frac{\al_1\ga_0-\al_0\ga_1-\dd_0\be_1+\be_0\dd_1}{r^2},\;
\dd=\frac{\al_1\dd_0-\al_0\dd_1-\ga_1\be_0+\ga_0\dd_1}{r^2},
\end{split}
\end{equation*} 
where $r^2=\al_1^2+\be_1^2+\ga_1^2+\dd_1^2$. If $(x_1,...,x_5)$ are the standard
coordinates on $\RR^5$, the (inverse) stereographic projection $\RR^4\rightarrow
S^4$ reads:  
\begin{equation}
(\al ,\be ,\ga , \dd)\mapsto (x_1,...,x_5)=\big(\frac{2\al}{r^2+1},
\frac{2\be}{r^2+1},
\frac{2\ga}{r^2+1}, \frac{2\dd}{r^2+1}, \frac{r^2-1}{r^2+1}\big). 
\end{equation}
Thus, if $x_3=x_4=0$ are the equations of a totally geodesic $S^2$ in $S^4$, 
the focal set of  $S^2$ in $S^4$ is the
$S^1$ given by the equations
$x_1=x_2=x_5=0$, corresponding to 
\begin{equation}\label{trei}\al=\be =0, \; \ga^2+\dd^2=1.\end{equation}

On the other hand, the moment map $\mu$ on  $\HP{1}$ can be
written in the 
non-homogeneous coordinate $h=h_0h_1^{-1}$ as $\mu=\ov{h}ih+i$.
Its level set $\mu^{-1}(0)$ is thus: 
\begin{equation}\label{patru}
\al_1^2+\be_1^2-\ga_1^2-\dd_1^2+1=0,\; \al\dd-\be\ga=0,\; \al\ga+\be\ga=0,
\end{equation}
and the systems \eqref{trei} and \eqref{patru} are equivalent.

We now treat the case $n>1$.  From the definition of the focal set, we only need to
look at geodesics normal to $\CP{n}$. These are normal to all the complex projective
lines $\CC L^1 \subset \CP{n}$, thus their focal points with respect to $\CP{n}$ are
also focal points for all the lines $\CC L^1$ contained in
$\CP{n}$. Any such line
$\CC L^1$ in
$\HP{n}$ belongs to a unique quaternionic projective line $\HH L^1$ and the latter
is totally geodesic in
$\HP{n}$. It follows that any geodesic that is normal to a $\CC L^1$ and tangent to
the corresponding
$\HH L^1$ at a given point,  remains tangent to $\HH L^1$ for its entire length.
Hence, if we show that through any point $x$ of $\CP{n}$ and for any $v\in T_x^\perp
\CP{n}$ there exists a projective quaternionic line $\HH L^1$ containing $v$, we may
deduce:
\bigskip
\par \noindent
{\it $Foc_{\HP{n}}\CP{n}$ is the union of all the focal sets $Foc_{\HH L^1}\CC L^1$.}
\bigskip
\par 
This can be seen from the diagram 
$$ \begin{array}{ccc}
 \CC^{n+1}\backslash \{0\} & \hookrightarrow & \HH^{n+1}\backslash \{0\}   \\
	\downarrow &    &  \downarrow  \\
\CP{n}& \hookrightarrow & \HP{n}. 
  \end{array}
 $$
by looking at the complex planes $L_2^\CC$ in
$\CC^{n+1}$ containing the fibre $\CC^*$ and the corresponding hypercomplex 2-planes 
$L_2^\HH$ in
$\HH^{n+1}$ containing the fibre $\HH^*$. Observe that, for any vector $\vv{v}\in
\HH^{n+1}$ which is normal
to the standard embedded $\CC^{n+1}$, there exist a  $L_2^\CC$ and a $L_2^\HH$ 
with $\text{span}_\RR\{L_2^\CC,\vv{v}\}\subset L_2^\HH$: indeed, if
$\HH^{n+1}=\text{span}_\HH \{\vv{e_0},...,\vv{e_n}\}$ and $\CC^{n+1}=\text{span}_\CC
\{\vv{e_0},...,\vv{e_n}\}$, we have:
$$\vv{v}\perp \CC^{n+1}\; \;  \text{if and only if}\; \;  \vv{v}=\sum_{a=0}^n
\la_aj\vv{e_a}+\mu_ak\vv{e_a}.$$ Hence, if
$$\vv{w}=-j\vv{v}=\sum_{a=0}^n\la_a\vv{e_a}-\mu_ai\vv{e_a},$$ we obtain 
$$L_2^\CC=\text{span}_\CC\{\vv{e_0}, \vv{w}\}$$ and $$
L_2^\HH=\text{span}_\RR\{\vv{e_0},i\vv{e_0}, j\vv{e_0}, k\vv{e_0}, \vv{w}, i\vv{w},
j\vv{w}, k\vv{w}\}$$ satisfy the condition. 

The identification of $Foc_{\HP{n}}\CP{n}$ with
$\mu^{-1}(0)$ is then completed by the following observation:
\bigskip
\par \noindent
{\it The subset $\mu^{-1}(0) \subset \HP{n}$ is the union of the zero sets of the
moment maps $\mu_1$ associated to the standard $U(1)$-action on all the projective
quaternionic lines  $\HH L^1 \subset
\HP{n}$.}
\bigskip
\par 
Let $\{p_0,p_1,...,p_n\}$ be the canonical frame of $\HP{n}$ with
unit point $u$ and let $[h_0:...:h_n]$ be the homogeneous coordinates with respect to
this frame. Accordingly, the moment map reads 
$\mu (h)=\sum_{a=0}^n\ov{h}_aih_a$. Fix a   $q\in \mu^{-1}(0)$, \emph{i.e.} 
$\sum_{a=0}^n\ov{q}_aiq_a=0$, and note that $q$ cannot be real; however, we
may suppose $q_0\neq 0$.  

Let $\HH L^1\cong \HP{1}$ be the quaternionic projective line through $q$ and $p_1$.
To compare the intersection $\mu^{-1}(0) \cap \HH L^1$ with the zero level set of the
moment map in $\HH L^1\cong \HP{1}$, we change the frame in 
$\HP{n}$  as follows. We want new homogeneous coordinates $[k_0:k_1:...:k_n]$ such
that
$q=[1:j:0:...:0]$ and the coordinates of $p_1$,..., $p_n$ remain unchanged. With the
new unit point $v=[q_0:1:-1:-1:...:-1]$, the
matrix of the change of coordinates is: 
$$A=\begin{pmatrix}q_0&0&0&\hdots&0\\
q_1-j&1&0&\hdots&0\\
q_2&0&-1&\hdots&0\\
\hdotsfor{5}\\
q_n&0&0&\hdots&-1
\end{pmatrix}.$$
Thus, $^t [h_0:...:h_n]= A ^t [k_0:...:k_n]$, and the
moment map is: 
\begin{equation*}
\mu=\sum_{a=0}^n\ov{h_a}ih_a
=\ov{k_0}(\sum_{a=0}^n\ov{q_a}iq_a)k_0-\ov{k_0}jijk_0-\sum_{b=1}^n\ov{k_b}ik_b
=-\sum_{b=0}^n\ov{k_b}ik_b.
\end{equation*}
It follows that $\mu^{-1}(0)$ is invariant under these projective changes
of coordinates. As the coordinates on $\HH L^1$
are $[k_0:k_1]$, its moment map is $\mu_1=-\ov{k_0}ik_0-\ov{k_1}ik_1$ with the same
level set
$\mu_1^{-1}(0)$ described for $\HP{1}$. This establishes the inclusion
$\mu^{-1}(0)Ê\vert_{\HH L^1} \subset \mu_1^{-1}(0)$. As the converse inclusion is
clear, for any projective line $\HH L^1$ in $\HP{n}$, the proof of the identification
$Foc_{\HP{n}}\CP{n}=\mu^{-1}(0)$ is complete. 

Next, we prove that $Foc_{\HP{n}}\CP{n}$ is isometric to the total space of the
induced Hopf bundle over the Grassmannian $Gr_2(\CC^{n+1})$. In fact, since
$Foc_{\HP{n}}\CP{n}$ is simply connected
(\cite{Ber 1}, p. 17), the existence of a diffeomorphism with the induced Hopf
bundle is a consequence of the following observation: 
\bigskip
\par \noindent
{\it Let $\pi:P\rightarrow Gr_2(\CC^{n+1})$ be a principal circle bundle with simply
connected $P$. Then $P$ is diffeomorphic to 
the total space of the induced Hopf bundle of $S^{2N+1}\rightarrow
\CC P^N$,
$N=\binom{n+1}{2}-1$, via the
Pl\"ucker embedding $Gr_2(\CC^{n+1}) \hookrightarrow  \CC P^N$.}
\bigskip
\par 
In fact, the principal $S^1$ -bundles over the base $B$ are classified by their
Chern class in
$H^2(B,\ZZ)$. Since
$H^2 (Gr_2(\CC^{n+1}),\ZZ)\cong \ZZ$ is generated by the class of the \ka\ form, one
can denote by
$P_a$ the $S^1$-bundle over  $Gr_2(\CC^{n+1})$ associated to $a\in \ZZ$. Observe that
$P_1$ (resp. $P_{-1}$) is the circle bundle associated to the canonical line bundle 
$O_{Gr_2(\CC^{n+1})}(1)$ (resp. its dual $O_{Gr_2(\CC^{n+1})}(-1)$). But 
dual complex line bundles are diffeomorphic  as real vector bundles. Hence 
$P_1$ is diffeomorphic to $P_{-1}$. 

Let us now show that if $P_a$ is simply connected, then $a=\pm 1$. From
$\pi_1(P)=0$, we
have   $H^1(P,\ZZ)=0$ and $H^2(P,\ZZ)$ torsion free. Thus the Gysin sequence 
of the $S^1$-bundle $\pi$:
\begin{equation*}
\begin{split}
0\rightarrow H^1(P,\ZZ)&\rightarrow H^0(Gr_2(\CC^{n+1}),\ZZ)\rightarrow
H^2(Gr_2(\CC^{n+1}),\ZZ)\\
&\rightarrow H^2(P,\ZZ)\rightarrow 0
\end{split}
\end{equation*}
reduces to:
$$0\longrightarrow \ZZ \stackrel{\cup c_1}{\longrightarrow} \ZZ
\longrightarrow H^2(P,\ZZ)\longrightarrow 0,$$
where $c_1$ is the Chern class of $\pi$. Hence,  $c_1\neq \pm 1$ implies  $
H^2(P,\ZZ)\cong
\ZZ_n$ for some $n\geq 2$. Consequently, $c_1$ is, up to sign, the Chern class of
the induced Hopf bundle. 
 
We now look at the metric $g_1$ inherited from 
$(S^{2N+1}, can)$ by the total space $Foc_{\HP{n}}(\CP{n})$ 
of the induced Hopf bundle and at the metric $g$ induced on
$Foc_{\HP{n}}(\CP{n})$ from $\HP{n}$. Note that  both $(Foc_{\HP{n}}(\CP{n}),g_1)$
and
$(Foc_{\HP{n}}(\CP{n}), g)$  are Riemannian submersions with geodesic fibres $S^1$
over
$Gr_2(\CC^{n+1})$ (cf. \cite{Or-Pi 3} for $g_1$ and \cite{Ber 1} for $g$). Thus 
Theorem 9.59 in \cite{Be}, p. 249, can be used to conclude that $g=g_1$. This
completes the proof of (i) and (ii) in  Theorem \ref{te1}.
\begin{re}
An alternative way of recognizing that $Foc_{\HP{n}}\CP{n}=\mu^{-1}(0)$ is to look
at the standard isometric action of
$SU(n+1)$ on $\HP{n}$ and to see that the two subsets are obtained as homogeneous
spaces of $SU(n+1)$ with the same isotropy groups. The homogeneous space is in fact
$\frac{SU(n+1)}{SU(2)\times SU(n-1)}$, which was shown in \cite{Ber 1}, p. 17, to be
a singular orbit of the action of $SU(n+1)$ (it was denoted there by
$Q^n$). A similar description  of $\mu^{-1}(0)$ as a homogeneous space is given in
\cite{Bat 1}, p. 65, in relation with the problem of studying local compatible
complex structures in $\HP{n}$. Compare also with 
\cite{Po}, p. 171, where $\frac{SU(n+1)}{SU(2)\times SU(n-1)}$ is called the
"Grassmannian of \emph{oriented} two-planes" of
$\CC^{n+1}$.
\end{re}
\begin{re}
We proved in \cite{Or-Pi 3} that the metric $g_1$ is Sasakian and $\eta$-Einstein,
and this is in accordance with formulas following Proposition 9 in \cite{Ber 1}. As
for the extrinsic geometry of
$Foc_{\HP{n}}(\CP{n})$, J. Berndt proves (Corollary 1 in \cite{Ber 1}) that the
immersion in $\HP{n}$ is minimal. It is interesting to observe that the same holds
for the immersion
$(Foc_{\HP{n}}\CP{n},g_1)$ in $(S^{2N+1}, can)$, see \cite{Or-Pi 3}. 
\end{re}
\subsection*{Proof of Theorem \ref{te2} (i)}   
	The identification of $\nu^{-1}(0)$ with the
intersection $M$ is an immediate consequence of the definition of $\nu$ and of the
first statement of Theorem \ref{te1}. 
	To see that $\nu^{-1}(0)$ is isometric to the
induced Hopf bundle discussed above, note that the last observation in the
proof of Theorem 3.1 (i) and (ii) still holds
for  principal circle bundles with simply connected $P$ over any complex algebraic
projective submanifold $B$ of
$\CP{N}$ with
$H^2(B,\ZZ)\cong \ZZ$. This applies in particular to $B=Z_{\widetilde{Gr}_4(\mathbb{R}^{n+1})}$, 
as soon as one recognizes that
$\nu^{-1}(0)$ is simply connected. To see this, observe that $\nu^{-1}(0)$ can be
regarded as the homogeneous space $\frac{SO(n+1)}{SO(n-3)\times Sp(1)}$ via the
transitive action of
$SO(n+1)$ on $\nu^{-1}(0)$. This last action comes, in fact, from the natural action
of
$SO(n+1)$ on $\HP{n}$ and   J. Berndt's observation that $SU(n+1)$ acts transitively
on
$Foc_{\HP{n}}(\CP{n})$. His argument can be triplicated to produce transitive actions
of  groups
$SU_i (n+1)$, $SU_j (n+1)$, $SU_k (n+1)$ (the first one is the standard $SU(n+1)$ in
$Sp(n+1)$, the other two are similarly defined by interchanging the r\^oles of the
unit quaternions $i$, $j$, $k$) on the zero level sets
$\mu_i^{-1}(0)$, $\mu_j^{-1}(0)$, $\mu_k^{-1}(0)$. The isotropy subgroup of the action
on $\nu^{-1}(0)$ at the point $[1:i:j:k:0...:0]$ is then $Sp(1)\times SO(n-3)$.
 The homotopy sequence associated to the 
homogeneous manifold $\frac{SO(n+1)}{Sp(1)\times SO(n-3)}$ thus shows that
$\nu^{-1}(0)$ is  simply connected, so that it is diffeomorphic with the total space
of the induced Hopf bundle, now over the Fano manifold
$Z_{\widetilde{Gr}_4(\RR^{n+1})}$.\\
As an  induced $S^1$-Hopf bundle, $\nu^{-1}(0)$ has a Sasakian
$\eta$ - Einstein structure. This can be
deformed to a Sasakian-Einstein metric, which is still a Riemannian submersion after
rescaling the standard K\"ahler-Einstein metric of
$Z_{\widetilde{Gr}_4(\mathbb{R}^{n+1})}$ \cite{Or-Pi 3}. On the other hand, the 
composition of the fiberings
$$\nu^{-1}(0)\stackrel{S^1}{\rightarrow} Z_{\widetilde{Gr}_4(\mathbb{R}^{n+1})}
\stackrel{S^2}{\rightarrow}\widetilde{Gr}_4(\mathbb{R}^{n+1})$$  is an
$SO(3)$-bundle which endows $\nu^{-1}(0)$ with a $3$-Sasakian structure
via the inversion Theorem 4.6 of \cite{Bo-Ga 1}.

\subsection*{Proof of Theorem \ref{te1} (iii) and Theorem \ref{te2} (ii)} 
In both cases  we construct a complex structure on the total 
space of an $S^3$ bundle over a Sasakian manifold, induced by the 
Hopf bundle $S^{4n+3}\rightarrow \HP{n}$. Both
Stiefel manifolds under discussion are induced Hopf $S^3$-bundles, as is easily
recognized by regarding them as homogeneous spaces. Thus the induced homogeneous
metrics make them Riemannian submersions with fibers $S^3$. More generally:
\begin{pr}
Let $\pi:P\rightarrow B$ be a principal $S^3$-bundle induced by the 
Hopf bundle $S^{4n+3}\rightarrow \HP{n}$, and let $g^B$, $\nabla^B$ be the induced metric 
and Levi Civita connection on $B \subset \HP{n}$. Assume that $B$ admits a Killing 
vector field $\xi$ such that $\varphi = \nabla^BÊ\xi$ defines on $B$ a Sasakian structure
Then $P$ admits an almost complex structure.
\end{pr}
\begin{proof}
Let $g^P$ be the induced metric on $P \subset S^{4n+3}$, so that $\pi$
is a Riemannian submersion. 
For any $X\in \mathcal{X}(B)$ we let $X^*$ be its horizontal lift on $P$.
Let $\xi_1$, $\xi_2$, $\xi_3$ be the unit Killing vector fields which give the 
usual $3$-Sasakian structure of the fibers $S^3$ and $\eta_1$, $\eta_2$, $\eta_3$ their 
duals with respect to the canonical metric of $S^3$. The $\xi_i$ may be viewed as
vector  fields on $P$. Let $\hat{\eta}_i$ be their dual forms with respect to the
metric
$g^P$.  Restricted to any fibre, the $\hat{\eta}_i$ coincide with the $\eta_i$. 
The usual 
splitting $TP=\mathcal{V}\oplus \mathcal{H}$ 
into vertical and horizontal parts may be refined to: 
$$TP=\text{span}\{\xi_1,\xi_2,\xi_3\}\oplus \text{span}\{\xi^*\}\oplus \mathcal{H}',$$
where $\mathcal{H}'$ represents the horizontal vector fields orthogonal to the horizontal lift 
$\xi^*$ of $\xi$..

Define the almost complex structure $J$ on $P$ by:
\begin{itemize}
\item $J\xi_1=\xi_2, \quad J\xi_2=-\xi_1$,
\item $J\xi_3=\xi^*, \quad J\xi^*=-\xi_3$,
\item $JX^*=(\varphi X)^*$ for any $X\in \mathcal{X}(B)$ orthogonal to $\xi$.
\end{itemize}
Note that for $X\perp \xi$, $X^*$ is a section of $\mathcal{H}'$. As
the restriction of $\varphi$ to $\xi^\perp$ is an endomorphism of
$\xi^\perp$, the last item in the definition is consistent.
One easily shows that $J^2=-1$ and $J$ is compatible with $g$. 
\end{proof}

To discuss the integrability of the constructed $J$, we follow the discussion developed 
\cite{Bo-Ga-Ma 3} for an almost hypercomplex structure, computing
the Nijenhuis tensor field on all possible combinations of vertical
and/or horizontal vector fields. 

Note first that the horizontal distribution $\mathcal{H}$ is an
$sp(1)$-connection in the induced Hopf $S^3$ bundle $P\rightarrow B$.
This follows from the fact that the bracket of any
horizontal $X^*$ with a vertical vector field is horizontal, a consequence of the Killing property 
of the $\xi_i$ with respect to $g^P$. Thus, in particular we get that for $i=1,2,3$:  
\begin{equation}\label{1}
 \hat{\eta}_k [\xi_i,X^*] = \hat{\eta}_k [\xi_i,\xi^*] =0,
\end{equation}
for $k=1,2,3$.
Now we can prove:
\begin{pr}
Assume that the
curvature form $\Omega$ of the $sp(1)$-connection $\mathcal{H}$
satisfies the following conditions:

1) $\Omega((\varphi X)^*,(\varphi Y)^*)=\Omega(X^*,Y^*)$,
	\emph{i.e.} $\Omega$ is of type $(1,1)$
with respect to $J$,  

2) $\Omega(X^*,\xi^*)=0$ for any $X\perp \xi$.

Then the almost complex structure $J$ is integrable.
\end{pr}

\begin{re}\label{r}
Conditions 1), 2) in the former Proposition express a compatibility between the Sasakian
structure of the base and the bundle structure of $P$. 
Since the vertical components of $\Omega$ are
precisely the $d\e_i$, the two conditions give
corresponding equations for $d\e_i$. Moreover, 
the condition $d\e_i((\varphi X)^*,(\varphi Y)^*)=d\e_i(X^*,Y^*)$ is 
easily checked to be 
equivalent with $d\e_i((\varphi X)^*,Y^*)+d\e_i(X^*,(\varphi Y)^*)=0$. 
\end{re}
Now we can give the proof of the Proposition 4.2:
\begin{proof}We  compute the Nijenhuis tensor field of $J$:
$$[J,J](A_1,A_2)=[A_1,A_2]+J[JA_1,A_2]+J[A_1,JA_2]-[JA_1,JA_2]$$
for all possible pairs $(A_1,A_2)$, 
noting that, 
due to the tensorial character of $[J,J]$, 
when dealing with horizontal (resp. vertical) vector fields it 
is enough to work with basic ones (resp. with $\xi_i$, $i=1,2,3$).

{\bf 1.} Let first $A_1=X^*$, $A_2=Y^*$, $X,Y\perp\xi$. Denoting by $\hat\eta$ the dual of $\xi^*$ we get:
 
\begin{equation}\label{one}
[X^*,Y^*]=[X^*,Y^*]'+\hat\eta([X^*,Y^*])\xi^*+\textrm{vertical part}.
\end{equation}
where the $'$ denotes the $\mathcal{H}'$ part. By
 $\pi$-corelation, $[X^*,Y^*]=[X,Y]^*{}'$. Moreover, from $\e(X^*)=\e(Y^*)=0$ we get
$\e([X^*,Y^*])=-d\e(X^*,Y^*).$ The vertical part of $[X^*,Y^*]$ must be of the form
$\sum_{i=1}^3\e_i([X^*,Y^*])\xi_i$.  Hence:
\begin{equation}\label{3}
[X^*,Y^*]=[X^*,Y^*]'-d\e(X^*,Y^*)\xi^*-\sum_{i=1}^3 d\e(X^*,Y^*)\xi_i.
\end{equation}
By similar computations:
\begin{equation*}
\begin{split}
J[JX^*,Y^*]&=(\varphi[\varphi X,Y])^*{}'+d\e((\varphi X)^*,Y^*)\xi_3-\\
&-d\e_1((\varphi X)^*,Y^*)\xi_2+d\e_2((\varphi
X)^*,Y^*)\xi_1-d\e_3((\varphi X)^*,Y^*)\xi^*,\\
J[X^*,JY^*]&=(\varphi[ X,\varphi Y])^*{}'+d\e( X^*,(\varphi Y)^*)\xi_3-\\
&-d\e_1( X^*,(\varphi Y)^*)\xi_2+d\e_2( X^*,(\varphi Y)^*)\xi_1-d\e_3(
X^*,(\varphi Y)^*)\xi^*,\\
[JX^*,JY^*]&=[\varphi X,\varphi Y]^*{}'-d\e( (\varphi X)^*,(\varphi
Y)^*)\xi^*\\
&-\sum_{i=1}^3 d\e_i((\varphi X)^*,(\varphi Y)^*)\xi_i.
\end{split}
\end{equation*}
Hence, using the $(1,1)$ character of $d\e_i$ and Remark \ref{r}, we find
\begin{equation*}
[J,J](X^*,Y^*)=[\varphi X,\varphi Y]^*{}'
-\left\{d\e(X^*,Y^*)-d\e((\varphi X)^*,(\varphi Y)^*)\right\}\xi^*.
\end{equation*}
As we know $[\varphi X,\varphi Y]+2d\eta(X,Y)\xi=0$ (this is the normality
condition of the Sasakian structure of $B$) the horizontal lift of
this (null) tensor field is zero, hence also its component in
$\mathcal{H}'$ is zero. But this is 
precisely $[\varphi X,\varphi Y]^*{}'$. 

On the other hand, on any Sasakian manifold one has:
$$d\eta(X,Y)=g(X,\varphi Y), \quad \varphi^2X=-X+\eta(X)\xi,$$
thus $d\eta(X,\varphi Y)+d\eta(\varphi X, Y)=0$ and $d\eta(\varphi X,\varphi
Y)-d\eta( X, Y)=0$. By horizontally lifting these equations  we
obtain the annulation of the $\xi^*$ component, hence $[J,J](X^*,Y^*)=0$. 

{\bf 2.} Let now  $A_1=X^*$, $A_2=\xi^*$ ($X\perp \xi$). Then: 
\begin{equation*}
\begin{split}
[J,J](X^*,
\xi^*)&=[X^*,\xi^*]+J[JX^*,\xi^*]+J[X^*,J\xi^*]-[JX^*,J\xi^*]=\\
&=[X^*,\xi^*]+J[(\varphi X)^*,\xi^*]-J[X^*,\xi_3]+[(\varphi
X)^*,\xi_3]=\\
&=[X^*,\xi^*]+J[(\varphi X)^*,\xi^*],
\end{split}
\end{equation*}
as the last two brackets are zero by \eqref{1}.
As above, using  Remark \ref{r} and $d\e_i(X^*,\xi^*)=0$ (condition 2)
in the statement) we obtain:
\begin{equation*}
[J,J](X^*, \xi^*)=([X,\xi]+\varphi[\varphi X,\xi])^*{}'=0,
\end{equation*}
because, as on a Sasakian manifold $\varphi \xi=0$,  we can
add in the  paranthesis the terms $[X,\varphi \xi]-[\varphi X,\varphi \xi]$
obtaining $([X,\xi]+\varphi[\varphi X,\xi]+\varphi[X,\varphi \xi]-[\varphi X,\varphi
\xi])^*{}'$ = $([\varphi,\varphi](X,\xi))^*{}'=0$ by the normality condition
on $B$.

{\bf 3.} We now choose $A_1=X^*$ and $A_2=\xi_i$ ($i=1,2$). We
have:
\begin{equation*}
[J,J](X^*,
\xi_1)=[X^*,\xi_1]+J[JX^*,\xi_1]+J[X^*,J\xi_1]-[JX^*,J\xi_1]=0,
\end{equation*}
because $\mathcal{H}$ is a $sp(1)$-connection. Similarly for $[J,J](X^*,
\xi_2)=0.$

{\bf 4.} For $A_1=X^*$ and $A_2=\xi_3$ we find:
\begin{equation*}
\begin{split}
[J,J](X^*,
\xi_3)&=[X^*,\xi_3]+J[JX^*,\xi_3]+J[X^*,J\xi_3]-[JX^*,J\xi_3]=\\
&=J[X^*,\xi^*]-[(\varphi X)^*,\xi^*],
\end{split}
\end{equation*}
 the brackets with $\xi_3$ being zero by \eqref{1}.  The horizontal component of the
remaining two brackets
 is $([\varphi[X,\xi]-[\varphi X,\xi])^*{}'-d\e_3(X^*,\xi^*)\xi_i$. The
$\xi^*$ component as well as the vertical component vanish by
assumption 2) in the statement. Finally, using Sasakian identities and
the normality condition on $B$ we
have $([\varphi[X,\xi]-[\varphi X,\xi])^*{}'=0$.   

{\bf 5.} Immediate computation shows that in the  remaining "mixed" case
$[J,J](\xi_i,\xi^*)=0$.

{\bf 6.} We are left with the computation of $[J,J]$ on
vertical fields. Obviously $[J,J](\xi_1,\xi_2)=0$. Then
\begin{equation*}
\begin{split}
[J,J](\xi_1,\xi_3)&=[\xi_1,\xi_3]+J[J\xi_1,\xi_3]+J[\xi_1,J\xi_3]-
[J\xi_1,J\xi_3]\\
&=[\xi_1,\xi3]+J[\xi_2,\xi_3]+J[\xi_1,\xi^*]-[\xi_2,\xi^*]=0
\end{split}
\end{equation*}
by $[\xi,\xi_j]=2\varepsilon_{ijk}\xi_k$ and \eqref{1}. Same 
arguments show that $[J,J](\xi_2,\xi_3)=0$, thus completing the proof.
\end{proof}
\begin{re}	
The K\"ahler form of $(P,g,J)$ is 
$$\omega=d\pi^*\eta +\pi^*\eta\wedge \eta_3-d\eta_3;$$
this shows that $d\omega\neq 0$, hence the structure is not K\"ahlerian.
A similar computation proves that $L_{\xi^*}J=L_{\xi_3}J=0$, thus $\xi^*$ and
$\xi_3$ are infinitesimal automorphisms of the constructed complex structure. 
\end{re}
\begin{re}
The complex structure $J$ on $P$ depends on the choice of a $3$-Sasakian 
structure on the fibre. But one can see that different choices of $3$-Sasakian 
triples $\{\xi_1,\xi_2,\xi_3\}$ produce complex structures that are conjugated 
in $End(TP)$.
\end{re}

We can now go back to the Stiefel manifolds $V_2(\CC^{n+1})$ and $\widetilde V_4(\RR^{n+1})$, and complete the proof of
Theorems 3.1 (iii) and 3.2 (ii). We have just to verify that the induced Hopf bundles 
$V_2(\CC^{n+1}) \rightarrow \mu^{-1}(0)$ and $\widetilde V_4(\RR^{n+1}) \rightarrow \nu^{-1}(0)$ inherit from the inclusions
$\nu^{-1}(0) \subset \mu^{-1}(0) \subset \HP{n}$ horizontal distributions $\mathcal H$ satisfying the curvature properties of
Proposition 4.2. Now property 1) simply express that the $sp(1)$-connection given by $\mathcal H$ is part of a
$u(2)$-connection in the bundles $V_2(\CC^{n+1}) \rightarrow Gr_2(\CC^{n+1})$ and $\widetilde V_4(\RR^{n+1}) \rightarrow 
Z_{\widetilde Gr_4(\RR^{n+1})}$, a fact easily recognized as in the case of hypercomplex structures in
$V_2(\CC^{n+1})$ (cf.\cite{Bo-Ga-Ma 2}, proof of Thm 1.10, as well as \cite{Bo-Ga-Ma 3}, Thm. 1.13). The meaning of property
2) is that the curvature of such a
$sp(1)$-connection is given by a (1,1)-form with respect to the almost complex structure
$J$. This follows for example from the computation carried out in \cite{Pi 1}, p. 63, for the Hopf bundle. One has to take
into account that the r\^ole of the $U(1)$ and of the $Sp(1)$-actions on both Stiefel manifolds correspond to the standard
basic choice of left and right multiplication by scalars on $\HH^{n+1}$. Then both on $V_2(\CC^{n+1})$ and on $\widetilde
V_4(\RR^{n+1})$ the almost complex structure $J$ satisfies the compatibility conditions with the Sasakian structures of  
$\mu^{-1}(0)$ and $\nu^{-1}(0)$, as expressed by conditions 1) and 2) of Proposition 4.2.
  Note that the complex structure obtained in this way on
$V_2(\CC^{n+1})$  projects to the complex structure 
underlying the K\"ahler structure of $Gr_2(\CC^{n+1})$.  
This latter is well known to be not compatible with the
quaternion K\"ahler structure of this Grassmannian. But it is precisely this
quaternion K\"ahler  structure which is lifted to a $3$-Sasakian structure and then,
by means of an appropriate circle bundle, to a hypercomplex structure on
$V_2(\CC^{n+1})$, cf. \cite{Bat 2}, \cite{Bo-Ga-Ma 2}. Hence our complex structure
is not  compatible with the standard hypercomplex one of $V_2(\CC^{n+1})$. 
This completes the proof of Theorem 3.1 (iii) and 3.2 (ii).

\subsection*{Proof of Corollary \ref{te3}}
It remains only to show that the total space of an induced $S^1$ Hopf bundle is
minimal in
$(S^{2N-1}, can)$. In fact, more generally, in a commutative diagram
$$\begin{CD}
\ov{N}@>\ov{i}>>\ov{M}\\
@V\pi_NVV @VV\pi_MV\\
N@>i>> M
\end{CD}
$$
of immersions $i$, $\ov{i}$ and Riemannian submersions $\pi_N$, $\pi_M$ with
totally geodesic fibres, we see that  
$N$ is minimal in $M$ if and only if $\ov{N}$  is minimal in $\ov{M}$.
This follows by a direct
computation of the mean curvature vector fields of $i$ and $\ov{i}$ using the Gauss
formula of a submanifold and formula (9.25 a) in \cite{Be}. 
\bigskip

\section{Further Observations}

As mentioned in the introduction, the zero level set $\mu^{-1}(0) \subset \HP{2}$ is
diffeomorphic to
a sphere $S^5$ and the projection to the reduced manifold can be identified with the
Hopf fibration
$S^5Ê\rightarrow \CP{2}$. Going to the next case $n=3$, the projection $\mu^{-1}(0)
\rightarrow Gr_2(\CC^{4})$, now an induced Hopf fibration, can be described by
looking at the families of submanifolds in $Gr_2(\CC^{4})$ that are either
K\"ahler-Einstein or quaternion K\"ahler. The Grassmannian $Gr_2(\CC^{4})$ admits in
fact a natural complex K\"ahler structure, as well as two distinct quaternion
K\"ahler structures induced via the isomorphism of vector bundles
$TGr_2(\CC^{4}) \cong V \otimes V^\perp$ from the (almost) hypercomplex structure on
the tautological vector bundle $V$ or on its orthogonal complement $V^\perp$ (cf.
\cite{Mor}). The families of submanifolds we want to look at on $Gr_2(\CC^{4})$ are
described as follows (see \cite{Ber 1}, \cite{Ber 2}, \cite{Ma}, \cite{Pi 2}). 
\par There are two families $\F$, $\F '$ of complex
projective  planes ${\bf C}P^2$, ${\bf C}P^{2'}$,  a family $\F ''$ of products
$\CP{1} \times
\CP{1}$, and a family $\F '''$ of spheres $S^4$; they are given by: 
\begin{equation*}
\begin{split}
\F =& \{\text{planes contained in a $3$-space of $\CC^4 $}\},\\
\F ' =& \{\text{planes through a
line of $\CC^4$} \},\\
\F '' =& \{\text{planes that are invariant for a  hypercomplex}\\ 
&\quad \text{structure $J$ of $\CC^4$}\},\\ 
\F ''' =& \{\text{planes given by pairs
of lines in  two fixed orthogonal}\\
& \quad \text{planes of ${\CC}^4$} \}.
\end{split}
\end{equation*}  

All these $4$-dimensional submanifolds of $Gr_2(\CC^{4})$ have nice intersection
properties,
some of which can also be formulated in terms of projective geometry of the lines in
the 3-dimensional space $\CP{3}$, the context where the Klein quadric $\CC Q^4\subset
\CP{5}$, isometric to $Gr_2(\CC^4)$, was first introduced. Instead of listing these
intersection properties on
$Gr_2(\CC^4)$ (cf. \cite{Ma}, pp. 508-512, for some of them; the remaining ones can
be similarly deduced), we formulate the corresponding properties for the families of
5-dimensional submanifolds obtained at the Sasakian Einstein level as induced Hopf
bundles over the members of families 
$\F$,..., $\F '''$.  
\begin{pr} The $9$-dimensional Sasakian Einstein focal set\\ 
$Foc_{\HP{3}}\CP{3}\cong\mu^{-1}(0) \subset \HP{3}$ contains the following families
of 5-dimensional submanifolds, each of which fibers in circles over a complex or a
quaternionic submanifold of
$Gr_2(\CC^{4})$. There are two families $\E$, $\E '$ of Sasakian 5-spheres $S^5$,
${S'}^{5}$, a family $\E ''$ of Sasakian products $S^3 \times S^2$ and a family $\E
'''$ of products $S^4 \times S^1$. The induced metrics on members of the families
$\E$, $\E '$,$\E ''$ are Sasakian $\eta$-Einstein, and can be modified to Sasakian
Einstein metrics by formula (2.1). The intersection properties of these 5-dimentional
submanifolds are the following.  
\par (a) Pairs of 5-spheres in the same family intersect in a circle, and pairs of
5-spheres in different families either do not intersect or intersect in an $S^3$. 
\par (b) Pairs of submanifolds of type $S^3 \times S^2$ intersect either in an $S^3$
or in a pair of disjoint circles. A submanifold of type $S^3 \times S^2$ intersects
a 5-sphere in an $S^3$.  \par (c) Any submanifold of type $S^4 \times S^1$
intersects a $5$-sphere in a circle, intersects an $S^3
\times S^2$ in two disjoint circles, and any pair of submanifolds of type $S^4
\times S^1$ intersect in two disjoint circles. 
\end{pr} 

This kind of geometry of submanifolds, now described for the level set  $\mu^{-1}(0) \subset \HP{3}$,
can be formulated for all the zero level sets of moment maps appearing in the
diagram of Corollary 3.1. There is such a level set for each odd dimension. For
example, $\nu^{-1}(0)
\subset {\HH}P^4$ is diffeomorphic to a sphere $S^7$, yielding as reduced manifold
$\HP{1}$ and fibering in circles over its twistor space ${\CC}P^3$, a K\"ahler
submanifold of the Grassmannian
$Gr_2({\CC}^5)$. Thus Sasakian ($\eta$)-Einstein submanifolds of type $S^5$ and $S^3
\times S^2$ can be determined in $S^7$ (cf. \cite{Or-Pi 3}), and intersection
properties like in Proposition 5.1 are obtained. 
\par The 11-dimensional example is the level set $\nu^{-1}(0) \subset {\HH}P^5$,
identified in Theorem 3.2 with an intersection of three focal sets in ${\HH}P^5$.
This 11-dimensional manifold is diffeomorphic to the (unique) 3-Sasakian homogeneous
manifold projecting in $SO(3)$ over the Wolf space $\widetilde Gr_4({\RR}^6)$, and
this latter manifold is isometric to
$Gr_2({\CC}^4)$. Thus again the geometry of 3-Sasakian and of Sasakian
($\eta$)-Einstein submanifolds
$\nu^{-1}(0) \subset {\HH}P^4$ is obtained from the same starting point as
Proposition 5.1.
\bigskip
\par Besides this kind of geometry of submanifolds, the Sasakian-Einstein level sets
$\mu^{-1}(0)$ share with some of the 3-Sasakian level sets $\nu^{-1}(0)$ a common
expression of their Poincar\'e polynomials. We have in fact:
\begin{pr} The Poincar\'e polynomial of $\mu^{-1}(0) \subset \HP{n}$ is given by:
$$Poin_{\mu^{-1}(0)}(t) = \sum_{i=0}^{[\frac{n-1}{2}]} (t^{4i}+t^{4n-3-4i})$$ 
\end{pr}
This is obtained by the Gysin sequence of the $S^1$-bundle $\mu^{-1}(0)
\rightarrow Gr_2({\CC}^{n+1})$, where the connecting homomorphism 
$$H^p(Gr_2({\CC}^{n+1})) \rightarrow H^{p+2}(Gr_2({\CC}^{n+1}))$$ is given by the
wedge product with the K\"ahler form of the Grassmannian. Since this wedge product
is injective up to the middle real dimension $p+2=2n-2$ (cf. \cite{Ga-Sa}, lemma
3.1, for the similar 3-Sasakian
situation), the Gysin sequence reduces to a series of short
exact sequences finishing with $$H^{2n-4}(Gr_2({\CC}^{n+1})) \rightarrow
H^{2n-2}(Gr_2({\CC}^{n+1})) \rightarrow H^{2n-2}(\mu^{-1}(0)).$$
This allows to compute the Betti numbers of $(\mu^{-1}(0))$ by differences of
consecutive even Betti numbers in $Gr_2({\CC}^{n+1})$. The Poincar\'e polynomial of
$Gr_2({\CC}^{n+1})$ is well known (see, for example, \cite{Bo-Tu}, p. 292) and by
writing it as:
\begin{equation} 
Poin_{Gr_2({\CC}^{n+1})}=(1+t^{2}+...+t^{2n-2})(1+t^{4}+t^{4m-4}), ~~~~(n+1=2m)
\end{equation}
\begin{equation}
Poin_{Gr_2({\CC}^{n+1})}=(1+t^{4}+...+t^{4m-4})(1+t^{2}+t^{2n}), ~~~~(n+1=2m+1)
\end{equation}
the conclusion is easily obtained.  
\par The following table gives the Poincar\'e polynomial of $\mu^{-1}(0) \subset
\HP{n}$ for low values of $n$ 
\begin{table}[!htb]  
{\small
\begin{center}
\begin{tabular}{|c|c|} 
\hline
 & \\Ê 
$n=3$ & $1+t^{4}+t^{5}+t^{9}$\\
& \\
$n=4$ & $1+t^{4}+t^{9}+t^{13}$\\     
 &  \\ 
$n=5$ & $1+t^{4}+t^{8}+t^{9}+t^{13}+t^{17}$ \\ 
& \\
$n=6$ & $1+t^{4}+t^{8}+t^{13}+t^{17}+t^{21}$\\     
 & \\
$n=7$ & $1+t^{4}+t^{8}+t^{12}+t^{13}+t^{17}+t^{21}+t^{25}$\\     
 & \\
$n=8$ & $1+t^{4}+t^{8}+t^{12}+t^{17}+t^{21}+t^{25}+t^{29}$\\     
 & \\
\hline
\end{tabular}
\end{center}
}

\end{table} 

The Poincar\'e polynomial of $\mu^{-1}(0)$ can be
compared with that of $\nu^{-1}(0)$, computed as for the homogeneous
3-Sasakian manifold $SO(n+1)/(SO(n-3) \times Sp(1))$. The latter has two different
expressions, according to whether $n+1$ is even or odd (see \cite{Ga-Sa} or
\cite{Bo-Ga 2}, p. 28). For odd values of $n+1=2k+3$ this expression is:
$$Poin_{\nu^{-1}(0)}(t) = \sum_{i=0}^{k-1} (t^{4i}+t^{8k-1-4i}),$$
and taking account of the dimensions, this is the same formula given in Proposition
5.2 for $Poin_{\mu^{-1}(0)}(t)$.       
\bigskip\null

\end{document}